\documentclass[11pt]{amsart}

\usepackage{amsmath,amsthm}
\usepackage {latexsym}
\usepackage{amssymb}

\newcommand{\les}{\lesssim}

\newcommand{\beeq}{\begin{equation}}
\newcommand{\eneq}{\end{equation}}

\newcommand{\bear}{\begin{eqnarray}}
\newcommand{\eear}{\end{eqnarray}}
\newcommand{\beq}{\begin{equation}}
\newcommand{\eeq}{\end{equation}}

\newcommand{\supp}{\mbox{\rm supp}}

\newcommand{\half}{\frac{1}{2}}

\newcommand{\eps}{{\varepsilon}}
\newcommand{\R}{{\mathbb R}}

\newcommand{\Compl}{{\mathbb C}}

\newcommand{\calg}{\,{\mathfrak g}}

\newcommand{\calW}{{\mathcal W}}

\newcommand{\calH}{{\mathcal H}}

\newcommand{\Laplace}{\triangle}

\newcommand{\kato}{{\mathcal K}}

\newcommand{\la}{\langle}
\newcommand{\ra}{\rangle}

\newcommand{\bdd}{{\mathcal B}}

\def\nn{\nonumber}
\def\calge1{\calg_{\vec{e_1}}}

\def\bm{\left( \begin{array}{cc}}
\def\endm{\end{array}\right)}
\def\ker{{\rm ker}}

\newtheorem{theorem}{Theorem}[section]
\newtheorem{lemma}[theorem]{Lemma}

\theoremstyle{remark}

\def\y{{y}}
\def\x{{x}}

\def\Bp{B^+}
\def\Btp{\tilde{B}^+}

\def\norm[#1][#2]{\Vert #1 \Vert_{#2}}

\def\Util3{\tilde{U}^{(3)}}

\def\Btp{\tilde{B}^+}
\def\Pac{P_{a.c.}}
\def\calK{{\mathcal K}}
\def\calU{{\mathcal U}}

\begin{document}

\title{Dispersive estimates for Schr\"odinger operators: A survey}
\author{W.\ Schlag}
\address{253-37 Caltech, Pasadena, CA 91125, U.S.A.}
\email{schlag@caltech.edu}
\thanks{The author was partially supported by the NSF grant
DMS-0300081 and a Sloan Fellowship. This article is based in part
on a talk that the author gave at the PDE meeting at the IAS in
Princeton in March of 2004. The author is grateful to the
organizers for the invitation to speak at that conference, as well
as to the Clay foundation and the IAS for their support. Also, he
wishes to thank Fabrice Planchon for useful comments on a
preliminary version of this article.}
\date{}

\maketitle

\section{Introduction}

The purpose of this note is to give a survey of some
recent work on dispersive estimates for the Schr\"odinger flow
\beeq
\label{eq:schr}
 e^{itH}P_{c}, \qquad H=-\Laplace+V \text{\ \ on\ \ }\R^d,\; d\ge1
\eneq
where $P_c$ is the projection onto the continuous spectrum of $H$.
$V$ is a real-valued potential that is assumed to satisfy some decay condition at infinity.
This decay is typically expressed in terms of the point-wise decay $|V(x)|\le C\la x\ra^{-\beta}$, for all $x\in\R^d$ and
some $\beta>0$. Throughout this paper, $\la x\ra=(1+|x|^2)^{\frac12}$.
Occasionally, we will use an integrability condition $V\in L^p(\R^d)$ (or a weighted variant thereof)
instead of a point-wise condition. These decay conditions will also be such that $H$ is asymptotically complete,
i.e.,
\[ L^2(\R^d)= L^2_{p.p.}(\R^d)\oplus L^2_{a.c.}(\R^d)\]
where the spaces on the right-hand side refer to the span of all eigenfunctions, and the absolutely continuous
subspace, respectively.

The dispersive estimate for \eqref{eq:schr} which we will be most concerned with is of the form
\beeq
\label{eq:disp}
\sup_{t\ne0}|t|^{\frac{d}{2}}\Big\|e^{itH}P_{c}f\Big\|_{\infty} \le C\|f\|_1 \text{\ \ for all \ \ }f\in L^1(\R^d)\cap L^2(\R^d).
\eneq
Interpolating with the $L^2$ bound $\Big\|e^{itH}P_{c}f\Big\|_{2}\le C\|f\|_2$ leads to
\beeq
\label{eq:disp'}
\sup_{t\ne0}|t|^{d(\frac12-\frac{1}{p})}\Big\|e^{itH}P_{c}f\Big\|_{p'} \le C\|f\|_p \text{\ \ for all \ \ }f\in L^1(\R^d)\cap L^2(\R^d).
\eneq
where $1\le p\le 2$.
It is well-known that via a $T^*T$ argument \eqref{eq:disp'} gives rise to the class of Strichartz estimates
\beeq
\label{eq:strich}
\Big\| e^{itH}P_{c}f\Big\|_{L^q_t(L^p_x)} \le C\|f\|_2, \text{\ \ for all \ \ } \frac{2}{q}+\frac{d}{p}=\frac{d}{2}, \; 2<q\le \infty.
\eneq
The endpoint $q=2$ is not captured by this approach, see Keel and Tao~\cite{KT}.

In heuristic terms, for the free problem $V=0$
the rate of decay $|t|^{-\frac{d}{2}}$ in \eqref{eq:disp} follows from $L^2$-conservation and the classical Newton law $\ddot{x}=0$
which leads to the trajectories $x(t)=vt+x_0$. Mathematically, \eqref{eq:disp} follows from the explicit solution
\[ (e^{-it\Laplace}f)(x) = C_d\,t^{-\frac{d}{2}} \int_{\R^d} e^{-i\frac{|x-y|^2}{4t}} \, f(y)\, dy.\]
For general $V\ne0$ no explicit solutions are available, and one needs to proceed differently.

If $V$ is small and $d\ge3$,
then one can proceed perturbatively. We will give examples of such arguments in Section~\ref{sec:3d}.
A purely perturbative approach cannot work in the presence of bound states of $H$ since
those need to be removed. In other words, in the presence of bound states the nature of the spectral measure
and/or resolvents of $H$ becomes essential.
 Since it is well-known that
bound states can arise for arbitrarily small potentials in dimensions $d=1,2$, see Theorem~XIII.11 in Reed and Simon~\cite{RS4},
we conclude that a perturbative approach will necessarily fail in those dimensions. On the other hand, if $d=3$
and $V$ satisfies the Rollnik condition
\[ \|V\|_{Roll}^2:=\int_{\R^6} \frac{|V(x)||V(y)|}{|x-y|^2}\, dxdy <\infty,\]
then Kato~\cite{kato} showed that $-\Laplace +  V$ is unitarily equivalent with $-\Laplace$ provided $4\pi\|V\|_{Roll}<1$.
Similar conditions are known for unitary equivalence if $d\ge4$.

Dispersive estimates for large $V$ and $d=3$ were established by Rauch~\cite{Rau} and Jensen, Kato~\cite{JK}.
In contrast to~\eqref{eq:disp}, these authors measured the decay on weighted $L^2(\R^3)$, i.e., they proved
that
\beeq
\label{eq:wdecay} \Big\| w e^{itH}P_{c} w f\Big\|_2 \le C|t|^{-\frac32}\|f\|_2
\eneq
with $w(x)=e^{-\rho\la x\ra}$ with some $\rho>0$ and $V$ exponentially decaying (Rauch) or $w(x)=\la x\ra^{-\sigma}$ for some $\sigma>0$ and
 $V$ decaying at a
power rate (Jensen, Kato). In addition, they needed to assume that the resolvent of $H$ has the property that
\beeq
\label{eq:res0} \limsup_{\lambda\to0} \| w (H-(\lambda\pm i0))^{-1} w \|_{2\to2} <\infty.
\eneq
This condition is usually referred to as {\em zero energy being neither an eigenvalue nor a resonance}.
While it is clear what it means for zero to be an eigenvalue of $H$, the notion of a resonance depends on
the norms relative to which the resolvent is required to remain bounded at zero energy, see~\eqref{eq:res0}.
In the context of $L^2$ with power weights, which are most commonly used, one says that there is a resonance at zero iff
there exists a distributional solution $f$ of $Hf=0$ with the property that $f\not\in L^2(\R^3)$ but
such that $\la x\ra^{-\sigma} f \in L^2(\R^3)$ for all $\sigma>\half$. With this definition the following holds: {\em \eqref{eq:res0} is
valid for $w(x)=\la x\ra^{-\half-\eps}$ for any $\eps>0$ iff zero is neither an eigenvalue nor a resonance.}
The proof proceeds via the Fredholm alternative and the mapping properties of $(-\Laplace + (\lambda+i0))^{-1}$ on weighted $L^2(\R^3)$ spaces,
see Section~\ref{sec:3d}. The notion of a resonance arises also in other dimensions, and we will discuss the cases $d=1,2$
in the corresponding sections below. If $|V(x)|\le C\la x\ra^{-2-\eps}$ with $\eps>0$ arbitrary, and $d\ge5$,
then $H$ cannot have any resonances at zero energy.
This is due to the fact that under these assumptions $(-\Laplace)^{-1}V\:: L^2(\R^d)\to L^2(\R^d)$.

Rauch and Jensen, Kato went beyond~\eqref{eq:wdecay} by showing that if zero is an eigenvalue and/or a resonance, then~\eqref{eq:wdecay}
fails. In fact, they observed that if zero is a resonance but not an eigenvalue, then
\[ C^{-1}< \sup_{\|f\|_2=1} \sup_{t\ge 1} |t|^{\half} \Big\| e^{itH}P_{c} f\Big\|_2 < C <\infty.\]
Furthermore, this loss of decay can occur also if zero is an eigenvalue {\em even though $P_c$ is understood
to project away the corresponding eigenfunctions.} They obtained these results as
corollaries of asymptotic expansions of $e^{itH}$ as $t\to\infty$ on weighted $L^2$ spaces.

These asymptotic expansions are basically obtained as the Fourier transforms of asymptotic expansions of the
resolvents (or rather, the imaginary part of the resolvents) around zero energy. In {\em odd} dimensions the latter are
of the form, with $\Im z>0$,
\beeq
\label{eq:res_exp} (-\Laplace + V - z^2)^{-1} = z^{-2} A_{-2} + z^{-1} A_{-1} +  A_0 + z A_{1} + O(z^2) \text{\ \ as \ }z\to0
\eneq
where the $O$-term is understood in the operator norm on a suitable weighted $L^2$-space. These expansions can of course
be continued to higher order $z^m$, with the degree of the weights in $L^2$ needed to control the error $O(z^m)$ increasing with~$m$.
In addition, the decay of $V$ needs to increase with~$m$ as well.
The operator $-A_{-2}$ is the orthogonal projection onto the eigenspace of $H$, and $A_{-1}$ is  a finite rank operator related to
both the eigenspace and the resonance functions.
In odd dimensions, the free resolvent $(-\Laplace + z^2)^{-1}$ is analytic for all $z\ne0$ (and if $d\ge3$ for all $z\in\Compl$), whereas in even dimensions
the Riemann surface of the free resolvent is that of the logarithm. In practical terms, this means that~\eqref{eq:res_exp}
needs to include (inverse) powers of $\log z$ in even dimensions.

In \cite{J1} and \cite{J2}, Jensen derived analogous expansions
for the resolvent around zero energy (and thus for the evolution
as $t\to\infty$) in  dimensions $d\ge4$. Resolvent expansion at
thresholds for the  cases $d=1$ and $d=2$ were treated by
Boll\'{e}, Gesztesy, Wilk~\cite{BGW}, and Boll\'{e}, Gesztesy,
Danneels \cite{BGD1}, \cite{BGD2}. However, their approach
requires separate treatment of the cases $\int V\, dx=0$ and $\int
V\,dx\ne0$. Moreover, for $d=2$ only the latter case was worked
out. A unified  approach to resolvent expansions was recently
found by Jensen and Nenciu in~\cite{JenNen}. Their method can be
applied to all dimensions, but in~\cite{JenNen} the authors only
present $d=1,2$ in detail, because for those cases novel results
are obtained by their method. The method developed by Jensen and
Nenciu was applied by Erdogan and the author for $d=3$,
see~\cite{ES1}, \cite{ES2}, and by the author for $d=2$,
see~\cite{Sch}. A very general treatment of resolvent expansions
as in~\eqref{eq:res_exp} and of local $L^2$ decay estimates can be
found in Murata's paper~\cite{Mur}. It is general in the sense
that Murata states expansions in all dimensions, and covers the
case of elliptic operators as well. However, his method is
partially implicit in the sense that the coefficients of the
singular powers in~\eqref{eq:res_exp} depend on operators that are
solutions of certain equations, but those equations are not solved
explicitly.

The first authors to address~\eqref{eq:disp} were Journe\'e,
Soffer, and Sogge~\cite{JSS}. Under suitable decay and
regularity conditions on $V$, and under the assumption that zero is neither
an eigenvalue nor a resonance they proved~\eqref{eq:disp}
 for $d\ge3$. In addition, they conjectured that~\eqref{eq:strich}
 should hold for all $V$ such that  $|V(x)|\le C\la x\ra^{-2-\eps}$
with arbitrary $\eps>0$ and for which $-\Laplace + V$ has neither an
eigenvalue nor a resonance at zero energy.

The decay rate $\la x\ra^{-2-\eps}$, which corresponds to
$L^{\frac{d}{2}}(\R^d)$ integrability, plays a special role in
dispersive estimates in particular, and the spectral theory of
$-\Laplace + V$ in general. On the one hand, potentials that decay
more slowly than $|x|^{-2}$ at infinity can lead to operators with
infinitely many negative bound states. On the other hand,
in~\cite{BPST1} and~\cite{BPST2} Burq, Planchon, Stalker, and
Tahvildar-Zadeh obtain Strichartz estimates for
\[ i\partial_t u + \Laplace u -\frac{a}{|x|^2}u =0, \]
provided $a>-(d-2)^2/4$ and $d\ge2$, and they show that this
condition is also necessary.  Furthermore, for the case of the
wave equation, it is known that point-wise decay estimates fail in
the attractive case $a<0$, see the work of Planchon, Stalker, and
Tahvildar-Zadeh.

For $d=3$ the assumptions on $V$ in~\cite{JSS} are  $|V(x)|\le C\la x\ra^{-7-\eps}$,
$\hat{V}\in L^1(\R^3)$, and some small
amount of differentiability of $V$. These requirements were subsequently relaxed by
Yajima~\cite{Y1}, \cite{Y2}, and \cite{Y3},
who proved much more, namely the $L^p$ boundedness of the wave operators
for $1\le p\le\infty$. A different approach, which lead
to even weaker conditions on $V$ was found by Rodnianski and the author~\cite{RS}
(for small $V$), as well as by
 Goldberg and the author~\cite{GS1}  (for large $V$).
In addition, the aforementioned conjecture from~\cite{JSS} is proved in~\cite{RS} (for large $V$).

Finally, Goldberg~\cite{gold2} proved that \eqref{eq:disp} --- and
not just~\eqref{eq:strich} --- holds for all $V$ for which
$|V(x)|\le C\la x\ra^{-2-\eps}$ with arbitrary $\eps>0$ and for
which $-\Laplace + V$ has neither an eigenvalue nor a resonance at
zero energy. In fact, he only required a suitable $L^p$ condition,
see Section~\ref{sec:3d}. In contrast, trying to adapt~\cite{GS1}
to higher dimensions has lead Goldberg and Visan~\cite{GolVi} to
show that for $d\ge4$, \eqref{eq:disp} {\em fails} unless $V$ has
some amount of regularity, i.e., decay alone is insufficient
for~\eqref{eq:disp} to hold if $d\ge4$. More precisely, they
exhibit potentials $V\in C_{\rm comp}^{\frac{d-3}{2}-\eps}(\R^d)$
for which the dispersive $L^1(\R^d)\to L^\infty(\R^d)$ decay with
power $t^{-\frac{d}{2}}$ fails.

The first results for $d=1$ are due to Weder~\cite{Wed1},
\cite{Wed2}, \cite{Wed3}, see also Artbazar and Yajima~\cite{AY}.
These authors make use of the following explicit expression for
the resolvent.  If $\Im z>0$, then
\[ (-\partial_x^2 + V - z^2)^{-1}(x,y) = \frac{f_+(x,z)f_{-}(y,z)}{W(z)} \text{\ \ if\ \ }x>y\]
and symmetrically if $x<y$. Here $f_{\pm}$ are the {\em Jost solutions} defined as solutions of
\[ -f_{\pm}''(\cdot,z) + Vf_{\pm}(\cdot,z) = z^2 f_{\pm}(\cdot,z)\]
with the asymptotics
\begin{align*}
 f_{+}(x,z) &\sim e^{ixz}  \text{\ \ as\ \ }x\to\infty \\
 f_{-}(x,z) &\sim e^{-ixz}  \text{\ \ as\ \ }x\to-\infty\ ,
\end{align*}
and $W(z)=W[f_{+}(\cdot,z),f_{-}(\cdot,z)]$ is their Wronskian.
These Jost solutions are known to exist and have boundary values as $\Im z\to0+$ as
long as $V\in L^1(\R)$ (in particular,
this proves that the spectrum of $H$ is purely a.c.~on $(0,\infty)$ for such~$V$).
In order for these boundary values $f_{\pm}(\cdot,\lambda)$ to be continuous at
$\lambda=0$ one needs to require that
$\la x\ra V(x)\in L^1(\R)$. In that case we say that {\em zero energy is a resonance} iff $W(0)=0$.
Note that the free case $V=0$ has a resonance at zero energy, since then $f_{\pm}(\cdot,0)=1$.
This condition is equivalent to the existence of a bounded solution $f$ of $Hf=0$
(in particular, zero cannot be an
eigenvalue).

Using some standard properties of the Jost solutions, see~\cite{DT}, Goldberg and the author proved
that
\beeq
\label{eq:1dest}
 \|e^{itH}P_c f\|_{L^\infty(\R)} \le C|t|^{-\frac12}\|f\|_{L^1(\R)}
\eneq provided $\la x\ra V(x)\in L^1(\R)$ and provided zero is not
a resonance. Note that in terms of pointwise decay, this is in
agreement with the $\la x\ra^{-2}$ threshold mentioned above. If
zero is a resonance, then the same estimate holds for all $V$ such
that $\la x\ra^2 V(x)\in L^1(\R)$. In Section~\ref{sec:1d} below,
we present a variant of~\eqref{eq:1dest} with faster decay that
seems to be new. It states that under sufficient decay on $V$ and
{\em provided zero is not a resonance}, \beeq \label{eq:nores}
 \|\la x\ra^{-1}e^{itH}P_{c}f \|_{L^\infty(\R)} \le C\,t^{-\frac32} \|\la x\ra f\|_{L^1(\R)}
\eneq
for all $t>0$. This estimate was motivated by the work of Murata~\cite{Mur} and
Buslaev and Perelman~\cite{BP} where
such improved decay was obtained on $L^2(\R)$ and with weights of the form $\la x\ra^{3.5+\eps}$.
It combines dispersive decay and the rate of propagation for~$H$.
However, to the best of the author's knowledge, \eqref{eq:nores} has not appeared before and we therefore
include a complete proof in Section~\ref{sec:1d}. A version of~\eqref{eq:nores} for the evolution of linearized
nonlinear Schr\"odinger equations
was crucial to the recent work~\cite{KS} by Krieger and the author on stable manifolds for all supercritical NLS
in one dimension.

Generally speaking, there is a very important difference between
the one-dimensional dispersive bounds and those in other
dimensions that have been proved so far, namely with regard to the
constants. Indeed, in the one-dimensional case these constants
exhibit an explicit dependence on the potential via the Jost
solutions, which are solutions to a Volterra integral equation. On
the other hand, in higher dimensions one resorts to a Fredholm
alternative argument in order to invert the operator
$H-(\lambda^2\pm i 0)$. This indirect argument is traditionally
used to prove the so-called limiting absorption principle for the
resolvent, see Agmon~\cite{agmon} and~\eqref{eq:lim_ap} below. Any
constructive proof of such an estimate for the perturbed resolvent
would be most interesting, as it would allow for quantitative
constants in dispersive estimates. Such a result was achieved by
Rodnianski and Tao, see~\cite{RodTao} as well as their article in
this volume. More generally, their work deals with dispersive
estimates for the Schr\"odinger operator on $\R^n$ (or other
manifolds) with variable metrics and is thus closely related to
the subject matter of this article. Unfortunately, it is outside
the scope of this review to discuss this exciting field of
research. For example, see Bourgain~\cite{master3},
Doi~\cite{Doi}, Burq, Gerard, Tzvetkov~\cite{BGT} (as well as
other papers by these authors), Hassell, Tao, Wunsch~\cite{HTW1},
\cite{HTW2}, Smith, Sogge~\cite{SS}, and Staffilani,
Tataru~\cite{ST}.

For the wave equation with a potential, dispersive estimates have
also been developed in recent years, see Cuccagna~\cite{cuc},
Georgiev and Visciglia~\cite{GV},  Pierfelice~\cite{Per},
Planchon, Stalker, and Tahvildar-Zadeh~\cite{PST1}, \cite{PST2},
d'Ancona and Pierfelice~\cite{AP}, as well as Stalker and
Tahvildar-Zadeh~\cite{StaTa}. There is some overlap with the
results here, in particular with respect to certain bounds on the
resolvent, but we will restrict ourselves to the Schr\"odinger
equation. For Klein-Gordon, see Weder's work~\cite{Wed4}.

\medskip Much of the work in this paper has been motivated by nonlinear problems
(see e.g.\ Bourgain's  book~\cite{Bbook}, in particular
pages~17--27).
In recent years there has been
much interest in the asymptotic stability of standing waves of the focussing NLS
\beeq
\label{eq:NLS}
 i\partial_t\psi + \Laplace\psi + f(|\psi|^2)\psi =0.
\eneq
A ``standing wave'' here refers to a solution of the form
$\psi(t,x)=e^{i\alpha^2 t}\phi(x)$ where $\alpha\ne0$ and
\beeq
\label{eq:phi}
 \alpha^2\phi-\Laplace\phi = f(\phi^2)\phi,
\eneq
or any solution obtained from this one by applying the symmetries of the
NLS, namely Galilei, scaling, and modulation
(if the nonlinearity is critical, then there is one more symmetry by the name of pseudoconformal).
Most work has been devoted to the standing wave generated by the ground state,
i.e., a positive, decaying, solution of~\eqref{eq:phi}.
In fact, such a solution must be radial and decay exponentially.
Linearizing~\eqref{eq:NLS} around a standing wave
yields a system of Schr\"odinger equations with non-selfadjoint matrix operator
\beeq
\label{eq:calH} \calH=\left[ \begin{matrix} -\Laplace + \alpha^2 - U & -W \\
W & \Laplace - \alpha^2 + U
          \end{matrix} \right ]
\eneq
and exponentially decaying, real-valued potentials $U$, $W$.
In order to address the question of asymptotic stability of standing waves, one needs to
study the spectrum of $\calH$, as well as prove dispersive estimates for $e^{it\calH}$ restricted to the stable
subspace (which is defined as the range of a suitable Riesz projection).
In the following sections we will mostly report on work on the {\em scalar} case rather than the system case.
However, most of what is being said can be generalized to systems, see e.g.~\cite{Cuc2}, \cite{RSS1},
\cite{ES2}, \cite{Sch2}, \cite{KS}.
Although it may seem that the exponential decay of the potential in~\eqref{eq:calH} may simplify
matters greatly, this turns out not to be the case. In fact, the method from the paper~\cite{GS1},
which is concerned with weakening the decay assumptions on~$V$  in the scalar, three-dimensional
case, has lead to the resolution of some open questions about
matrix operators as in~\eqref{eq:calH}, see~\cite{ES2},
\cite{Sch2}, \cite{KS}.

\section{Dimensions three and higher}
\label{sec:3d}

We start with a perturbative argument for small $V$ that can be considered
as a sketch of the method from~\cite{JSS}.
As above, let $H=-\Laplace+V$ and suppose $d\ge3$. Define
\[ M_0 = \sup_{0\le t} \sup_{\|f\|_{1\cap 2}=1}\la t\ra^{\frac{d}{2}} \|e^{itH_0} f\|_{2+\infty}, \qquad M(T)=\sup_{0\le t\le T}\sup_{\|f\|_{1\cap 2}=1} \la t\ra^{\frac{d}{2}} \|e^{itH} f\|_{2+\infty}.\]
Here
\[ \|f\|_{1\cap 2} = \|f\|_{L^1\cap L^2},\qquad \|f\|_{2+\infty} = \inf_{f_1+f_2=f} (\|f_1\|_2+\|f_2\|_\infty).\]
Then the Duhamel formula
\[
e^{itH}  = e^{itH_0} + i \int_0^t e^{i(t-s)H_0} V e^{isH}\, ds
\]
implies that
\[
 M(T) \le M_0 + \la T\ra^{\frac{d}{2}}\int_0^T M_0\la t-s\ra^{-\frac{d}{2}} \|V\|_{1\cap \infty} M(T)\la s\ra^{-\frac{d}{2}}\, ds
 \le M_0 + C\,\|V\|_{1\cap \infty}M_0M(T).
\]
Consequently, as long as \[C\,\|V\|_{1\cap \infty}M_0\le \frac12,\]
we obtain the bound \[\sup_{T\ge0}M(T)\le 2M_0.\]
Note first that such an argument necessarily fails if $d=1,2$ due
to the non-integrability of $t^{-\frac{d}{2}}$ at infinity.
Moreover, the are spectral reasons for this failure which we outlined in the introduction.
Second, we would like to
point out that it equally applies to time-dependent potentials provided the
evolution $e^{itH}$ is replaced with the
propagator of the associated Schr\"odinger equation. The inclusion of the space $L^2$ allows us to deal with
the singularity of $t^{-\frac{d}{2}}$ at $t=0$ which arises in the $L^1\to L^\infty$
estimate. In order to avoid it, Journe\'{e}, Soffer, and
Sogge use the bound \[\|e^{-itH_0} Ve^{itH_0}\|_{p\to p}\le \|\hat{V}\|_1,\]
which holds uniformly in $1\le p\le \infty$. This explains the origin of the condition $\hat{V}\in L^1$ in their paper.

The main difficulty in~\cite{JSS} is of course the fact that $V$ is large.
Let us first present an unpublished
argument of Ginibre~\cite{Gin} in dimensions $d\ge3$
that allows passing from the weighted (or local) decay~\eqref{eq:wdecay}
to global decay, albeit in the form of a
$L^1\cap L^2 \to L^2+L^\infty$ estimate rather than the one in~\eqref{eq:disp}.
Applying the Duhamel formula twice, we obtain
\begin{align}
e^{itH} P_c &= e^{itH_0}P_c + i\int_0^t e^{i(t-s)H_0} V P_c e^{isH_0}\, ds \nn \\
& - \int_0^t\int_0^s e^{i(t-s)H_0} V e^{i(s-\sigma)H}P_c V e^{i\sigma H_0}\, d\sigma\, ds\ . \label{eq:2duh}
\end{align}
As long as $V$ decays sufficiently rapidly so as to absorb the weights $w$, i.e., such that
\[ \|w^{-1} V\|_{L^1\cap L^\infty} < \infty\ ,\]
we can combine the  $L^1\cap L^2 \to L^2+L^\infty$  bound
\beeq
\label{eq:free} \|e^{itH_0} \|_{2+\infty} \le C\la t\ra^{-\frac{d}{3}} \|f\|_{1\cap 2}
\eneq
with \eqref{eq:wdecay} as above to conclude from \eqref{eq:2duh} that \eqref{eq:free} also holds for $H$.
Here we also used that $P_c:L^1\cap L^\infty\to L^1\cap L^\infty$ which holds
provided all eigenfunctions of $H$ with negative
eigenvalue belong to $L^1\cap L^\infty$ (recall that zero is assumed not to be an eigenvalue).
That property, however, follows from Agmon's exponential decay bound~\cite{Agm2}
and Sobolev imbedding provided $V$ also
has some small of regularity.

This argument, however, does not shed much light on the question of $L^1\to L^\infty$
bounds (without assuming more regularity on $V$).
The inclusion of $L^2$ is undesirable for a number of reasons, the main one being nonlinear applications.
We therefore proceed differently, and first recall the small-potential argument from~\cite{RS} in $d=3$. For
certain standard details we refer the reader to~\cite{RS}.

The starting point is the standard fact
\beeq
\label{eq:H_rep} e^{itH}P_{ac} = \int_0^\infty e^{it\lambda} E_{ac}(d\lambda),
\eneq
where $E_{ac}$ is the absolutely continuous part of the spectral resolution. Its density is given by
\[ \frac{d E_{ac}(\lambda)}{d\lambda} = \frac{1}{2\pi i} [(H-(\lambda+i0))^{-1}  - (H-(\lambda-i0))^{-1}]  \]
on $\lambda>0$. As already mentioned before, Kato's theorem~\cite{kato} insures that $E_{ac}=E$ provided
\beeq
\label{eq:rollnik} \|V\|_{R}^2:=\int_{\R^3\times \R^3} \frac{|V(x)|\,|V(y)|}{|x-y|^2}\,dx\,dy < (4\pi)^2\ .
\eneq
Let $R_V(z)=(-\Laplace+V-z)^{-1}$ and $R_0(z)=(-\Laplace-z)^{-1}$.
Then with $V$ as in \eqref{eq:rollnik}, for all $f, g\in L^2(\R^3)$ and $\eps\ge0$
one has the Born series expansion
\begin{equation}
\label{eq:Resseries}
\langle R_V(\lambda\pm i\eps) f,\,g\rangle - \langle R_0(\lambda\pm i\eps) f,\,g\rangle =
\sum_{\ell=1}^\infty (-1)^\ell \langle R_0(\lambda\pm i\eps)(VR_0(\lambda\pm i\eps))^\ell f,g\rangle
\end{equation}

It is well-known that the resolvent $R_0(z)$ for $\Im z\ge0$ has the kernel
\begin{equation}
\label{eq:R0kernel}
R_0(z)(x,y) = \frac{\exp(i\sqrt{z}|x-y|)}{4\pi|x-y|}
\end{equation}
with  $\Im(\sqrt{z})\ge0$.
Then there is the following simple lemma that is basically an instance of stationary phase.
For the proof we refer the reader to~\cite{RS}.

\begin{lemma}
\label{lem:statphas}
Let $\psi$ be a smooth, even  bump function
with $\psi(\lambda)=1$ for $-1\le\lambda\le 1$ and $\supp(\psi)\subset[-2,2]$.
Then for all $t\ge1$ and any real~$a$,
\begin{equation}
\label{eq:decay}
\sup_{L\ge 1}\Bigl| \int_0^\infty e^{it\lambda} \sin(a\sqrt{\lambda})\,
\psi(\frac{\sqrt\lambda}{L})\,d\lambda\Bigr| \le
C \,t^{-\frac32}\,|a|
\end{equation}
where $C$ only depends on ~$\psi$.
\end{lemma}

In addition to \eqref{eq:rollnik}, we will assume that
\beeq
\label{eq:kato} \|V\|_{\kato}:=\sup_{x\in\R^3} \int_{\R^3} \frac{|V(y)|}{|x-y|}\;dy < 4\pi
\eneq
In \cite{RS} this norm was introduced by the name of {\em global Kato norm}
(it is closely related to the well-known {\em Kato norm},
see Aizenman and Simon~\cite{AS}, \cite{barry}).
The following lemma explains to some extent why condition~\eqref{eq:kato} is needed.
Iterated integrals as in~\eqref{eq:multint} will appear in a series expansion of the
spectral resolution of $H=-\Laplace+V$.
For the sake of completeness, and in order to show how these global
Kato norms arise, we reproduce the simple proof from~\cite{RS}.

\begin{lemma}
\label{lem:iter}
For any positive integer $k$ and $V$ as above,
\begin{equation}
\label{eq:multint}
\sup_{x_0,x_{k+1}\in\R^3}\int_{\R^{3k}}
\frac{\prod_{j=1}^k |V(x_j)|}{\prod_{j=0}^k|x_j-x_{j+1}|}\sum_{\ell=0}^k |x_\ell-x_{\ell+1}|
\; dx_1\ldots\,dx_k \le (k+1) \|V\|_{\kato}^k.
\end{equation}
\end{lemma}
\begin{proof}
Define the operator ${\mathcal A}$ by the formula
$$
{\mathcal A} f (x) = \int_{\R^3} \frac{|V(y)|}{|x-y|} \,f(y)\,dy.
$$
Observe that the assumption \eqref{eq:kato} on the potential $V$ implies
that ${\mathcal A}:\,L^\infty\to L^\infty$ and
$\|{\mathcal A}\|_{L^\infty\to L^\infty}\le c_0$ where we have set $c_0:=\|V\|_{\kato}$ for convenience.
Denote by $<,>$ the standard $L^2$ pairing.
In this notation the estimate \eqref{eq:multint} is equivalent to proving
that the operators ${\mathcal B}_k$ defined as
$$
{\mathcal B}_k f = \sum_{m=0}^k <f,{\mathcal A}^{k-m} 1>{\mathcal A}^m 1
$$
are bounded as operators from $L^1\to L^\infty$ with the bound
$$
\|{\mathcal B}_k\|_{L^1\to L^\infty}\le (k+1) c_0^k.
$$
For arbitrary $f\in L^1$ one has
$$
\aligned
\|{\mathcal B}_k f\|_{L^\infty}&\le  \sum_{m=0}^k |<f,{\mathcal A}^{k-m} 1>|\,\,
\|{\mathcal A}^m 1\|_{L^\infty}\\ & \le \sum_{m=0}^k \|{\mathcal A}^{k-m}\|_{L^\infty\to L^\infty}
 \|{\mathcal A}^{m}\|_{L^\infty\to L^\infty} \|f\|_{L^1} \\ &\le
\sum_{m=0}^k c_0^k \|f\|_{L^1}\le (k+1) c_0^k \|f\|_{L^1},
\endaligned
$$
as claimed.
\end{proof}

We are now in a position to prove the small $V$ result from~\cite{RS}. In~\cite{Per}, Perfelice obtained an
analogous result for the wave equation.

\begin{theorem}
\label{thm:high}
With $H=-\Laplace+V$ and~$V$
satisfying the conditions \eqref{eq:rollnik} and \eqref{eq:kato},  one has the bound
\[ \Bigl\|e^{itH}\Bigr\|_{L^1\to L^\infty} \le C\, t^{-\frac32}\]
in three dimensions.
\end{theorem}
\begin{proof}
Fix  a real potential~$V$ as above, as well as any $L\ge1$, and real $f,g\in C_0^\infty(\R^3)$.
Then applying \eqref{eq:Resseries}, \eqref{eq:R0kernel},
Lemma~\ref{lem:statphas}, and Lemma~\ref{lem:iter} in this order,
we obtain
\bear
&& \sup_{L\ge1} \Bigl|\bigl \langle e^{itH}\psi(\sqrt{H}/L) f,g \bigr\rangle \Bigr| \nonumber\\
&\le&
\sup_{L\ge1} \Bigl| \int_0^\infty e^{it\lambda}\;\psi(\sqrt{\lambda}/L) \langle E'(\lambda)f,g
\rangle \,d\lambda \Bigr| \nonumber \\
&=& \sup_{L\ge1} \Bigl|\int_0^\infty e^{it\lambda}\, \psi(\sqrt{\lambda}/L) \Im
\langle R_V(\lambda+i0)f,g \rangle \,d\lambda \Bigr| \nonumber\\
&=& \sup_{L\ge1} \Bigl| \int_0^\infty e^{it\lambda}\;\psi(\sqrt{\lambda}/L)
\sum_{k=0}^\infty \Im \langle R_0(\lambda+i0)(VR_0(\lambda+i0))^k\,f,g \rangle \,
d\lambda \Bigr| \nonumber
\eear
To proceed, we now use the explicit form of the free resolvent. This yields
\bear
&\le& \sum_{k=0}^\infty \int_{\R^6} |f(x_0)||g(x_{k+1})|
\int_{\R^{3k}} \frac{\prod_{j=1}^k |V(x_j)|}{\prod_{j=0}^k 4\pi |x_j-x_{j+1}|}\cdot\nonumber\\
&& \qquad\qquad\qquad \cdot \sup_{L\ge1} \Bigl| \int_0^\infty e^{it\lambda}\;
\psi(\sqrt{\lambda}/L) \sin\Bigl(\sqrt{\lambda}\sum_{\ell=0}^k |x_\ell-x_{\ell+1}|\Bigr)\,
d\lambda \Bigr|
\; d(x_1,\ldots,x_k)\,dx_0\,dx_{k+1} \label{eq:gross}\\
&\le& Ct^{-\frac32} \sum_{k=0}^\infty \int_{\R^6} |f(x_0)||g(x_{k+1})|
\int_{\R^{3k}} \frac{\prod_{j=1}^k |V(x_j)|}{(4\pi)^{k+1}\prod_{j=0}^k|x_j-x_{j+1}|}
\sum_{\ell=0}^k |x_\ell-x_{\ell+1}|
\; d(x_1,\ldots,x_k)\;dx_0\,dx_{k+1} \nonumber \\
&\le& Ct^{-\frac32} \sum_{k=0}^\infty \int_{\R^6} |f(x_0)||g(x_{k+1})|\;
(k+1) (\|V\|_{\kato}/4\pi)^k \;dx_0\,dx_{k+1}\nonumber \\
&\le& Ct^{-\frac32} \|f\|_1\|g\|_1\ , \nonumber
\eear
since $\|V\|_{\kato}<4\pi$. In order to pass to \eqref{eq:gross}
one uses the explicit representation
of the kernel of $R_0(\lambda+i0)$, see~\eqref{eq:R0kernel},
which leads to a $k$-fold integral. Next, one interchanges the
order of integration in this iterated integral.
\end{proof}

The next step is to remove the smallness assumption on~$V$.
This was done in~\cite{GS1} for potentials
decaying like $|V(x)|\le C\la x\ra^{-\beta}$ with $\beta>3$.
The proof required splitting the energies into
the regions $[\lambda_0,\infty)$ (the ``large energies'')
and $[0,\lambda_0]$ (the ``small energies'') where $\lambda_0>0$ is small.
In the regime of large energies, one expands the resolvent $R_V$ into a finite Born series
\begin{align}
R_V(\lambda^2\pm i0) &=  \sum_{\ell=0}^{2m+1} R_0(\lambda^2\pm i0)(-VR_0(\lambda^2\pm i0))^\ell \nn \\
& + R_0(\lambda^2\pm i0)(VR_0(\lambda^2\pm i0))^mV
R_V(\lambda^2\pm i0)V(R_0(\lambda^2\pm i0)V)^mR_0(\lambda^2\pm i0)
\label{eq:res_ident}
\end{align}
where $m$ is any positive integer. All but the last term
(which involves $R_V$) is treated by the same argument
from~\cite{RS} that we sketched previously.  To bound the
contribution of the final term in~\eqref{eq:res_ident},
let $R_0^{\pm}(\lambda^2):= R_0(\lambda^2\pm i0)$.
Moreover, set
\[ G_{\pm,x}(\lambda^2)(x_1):= e^{\mp i\lambda|x|}R_0(\lambda^2\pm i0)(x_1,x)
= \frac{e^{\pm i\lambda(|x_1-x|-|x|)}}{4\pi|x_1-x|}. \]
Similar kernels appear already in Yajima's work~\cite{Y4} (see his high
energy section). Hence, we are led to proving that
\begin{align}
& \left|
\int_0^\infty e^{it\lambda^2}e^{\pm i\lambda(|x|+|y|)}\; (1-\chi(\lambda/\lambda_0))
\lambda \Big\la VR^{\pm}_V(\lambda^2)V (R_0^{\pm}(\lambda^2)V)^m
G_{\pm,y}(\lambda^2), (R_0^{\mp}(\lambda^2)V)^m G_{\pm,x}^*(\lambda^2) \Big\ra \, d\lambda
\right| \nn \\
& \les |t|^{-\frac32} \label{eq:main}
\end{align}
uniformly in $x,y\in\R^3$. Here $\chi$ is a bump function which is
equal to one on a neighborhood of the origin.
The estimate~\eqref{eq:main} is proved by means of stationary phase and the
{\em limiting absorption principle}. The latter
refers to estimates of the form, with $\lambda>0$ and $\sigma>\frac12$,
\beeq
\label{eq:lim_ap}
 \|R_0(\lambda^2\pm i0) f\|_{L^{2,-\sigma}} \le C(\lambda)\, \|f\|_{L^{2,\sigma}},
\eneq
where $L^{2,\sigma}=\la x\ra^{-\sigma} L^2$, see~\cite{agmon}.
Similar estimates also hold for the derivatives of $R_0$ in $\lambda$.
Moreover, $C(\lambda)$ decays power-like with $\lambda\to\infty$. By means of
the resolvent identity and arguments of Agmon and Kato analogous
estimates hold for $R_V(\lambda^2\pm i0)$ (this essentially amounts to
the absence of imbedded eigenvalues in the continuous spectrum).
These properties insure that the  integrand in~\eqref{eq:main}, viz.
\[  a_{x,y}(\lambda):=(1-\chi(\lambda/\lambda_0))
\lambda \Big\la VR^{\pm}_V(\lambda^2)V (R_0^{\pm}(\lambda^2)V)^m
G_{\pm,y}(\lambda^2), (R_0^{\mp}(\lambda^2)V)^m G_{\pm,x}^*(\lambda^2) \Big\ra
\]
decays at least as fast as $\lambda^{-2}$ (provided $m$ is large) and
is twice differentiable, say.  Moreover, due to the presence of the
functions $G_{\pm,y}$ and $G_{\pm,x}^*$ at the edges, one checks that
if the critical point $\lambda_1=\frac{|x|+|y|}{2t}$ of the phase falls
into the support of this integrand, which requires $\lambda_1\ge\lambda_0$,
then the entire integral is bounded by
\[ t^{-\frac12}|a_{x,y}(\lambda_1)|\le Ct^{-\frac12}(\la x\ra \la y\ra)^{-1} \le Ct^{-\frac32},\]
as desired.

In the low-energy regime $\lambda\in[0,\lambda_0]$, one writes
\begin{align}
& \Big \la e^{itH} \chi(\sqrt{H}/\lambda_0)\,\Pac \, f,g \Big\ra
\nn \\
& =  \int_0^\infty
e^{it\lambda^2}\lambda\, \chi(\lambda/\lambda_0)\, \Big \la [R_V(\lambda^2+i0)-R_V(\lambda^2-i0)]f,g \Big\ra
\, \frac{d\lambda }{\pi i} \nn
\end{align}
and proceeds via the resolvent identity
\begin{equation}
\nn
R_V^\pm(\lambda^2) = R_0^\pm(\lambda^2) - R_0^\pm(\lambda^2)V
 (I + R_0^\pm(\lambda^2)V)^{-1} R_0^\pm(\lambda^2).
\end{equation}
Expanding $R_0^\pm(\lambda^2)$ around zero, the invertibility of  $I + R_0^\pm(\lambda^2)V$
reduces to the invertibility of
\beeq
\label{eq:S0}
S_0:=I + R_0^\pm(0)V.
\eneq
 However, the latter is equivalent to zero energy being neither an eigenvalue nor a resonance.
Writing $R_0^\pm(\lambda^2) = R_0(0) + B^\pm(\lambda)$, we conclude that
$$[I+R_0^\pm(\lambda^2)V]^{-1} = S_0^{-1}[I+B^\pm(\lambda)VS_0^{-1}]^{-1}=: S_0^{-1} \tilde{B}^{\pm}(\lambda).$$
Some elementary calculations based on the explicit form  of the kernel of $R_0$
and the decay of $V$ then reduce the $t^{-\frac32}$
dispersive decay to the finiteness of
\[
 \int_{-\infty}^\infty
 \norm[[\chi_0(\Btp)']^\vee(u)][HS(-1^-,-2^-)]\, du  \text{\ \ and\ \ }
 \int_{-\infty}^\infty
 \norm[[\chi_0 \Btp]^\vee(u)][HS(-1^-,-2^-)]\, du
\]
where the norm is that of the Hilbert-Schmidt operators from $L^{-1-\eps}(\R^3)\to L^{-2-\eps}(\R^3)$.
Expanding into a Neuman series
$$\Btp(\lambda) \ = \ [I + \Bp(\lambda)VS_0^{-1}]^{-1}\  =\  \sum_{n=0}^\infty
 \big(- \Bp(\lambda)VS_0^{-1}\big)^n$$
and making careful use of  the explicit kernel
\[
B^\pm(\lambda)(\x,\y) = \frac{e^{\pm i\lambda|\x-\y|} - 1}{4\pi |\x-\y|}
\]
finishes the proof, see~\cite{GS1}.

\bigskip This argument was extended in various directions. First Yajima~\cite{Y5} and independently,
Erdogan and the
author~\cite{ES1}, have adapted it to the case of zero energy being an eigenvalue and/or a resonance.
The difference is of course that in this case $S_0$ as in~\eqref{eq:S0} is no longer invertible
and $(I + R_0^\pm(\lambda^2)V)^{-1}$ involves singular powers of $\lambda$. Yajima uses the
expansion from~\cite{JK} for
this purpose, whereas~\cite{ES1} use the method from~\cite{JenNen}.
The latter is based on the symmetric resolvent
identity and is therefore entirely situated in $L^2$ rather than weighted $L^2$.
The following theorem is from~\cite{ES1}. Yajima proves the same, but assuming less decay on $V$.

\begin{theorem}\label{T:scalar1}
Assume that $V$ satisfies $|V(x)|\leq C \langle x\rangle^{-\beta}$ with $\beta>10$
and assume that there is a resonance at energy zero
but that zero is not an eigenvalue. Then there is a time dependent rank one operator $F_t$ such that
$$
\left\|e^{itH} P_{ac}-t^{-1/2} F_t  \right\|_{1\rightarrow\infty}\leq C  t^{-3/2},
$$
for all $t>0$ and $F_t$ satisfies
\[
\sup_t\left\|F_t\right\|_{L^1\rightarrow L^\infty}<\infty, \qquad
\limsup_{t\to\infty}\left\|F_t\right\|_{L^1\rightarrow L^\infty} >0.
\]
A similar result holds also in the presence of eigenvalues, but in general $F_t$ is no longer of rank one.
\end{theorem}

The paper~\cite{ES2} extends these methods further, namely to the case of systems of the type that arise
from linearizing NLS around a ground state standing wave.

In another direction, Goldberg has improved on the method from~\cite{GS1} in several aspects.
In~\cite{gold1}, he proves that $V\in L^{\frac32(1+\eps)}(\R^3)\cap L^1(\R^3)$ suffices for
the dispersive estimate
(assuming of course that zero is neither and eigenvalue nor a resonance).
This amount of integrability is analogous to the $\beta>3$ point-wise decay from~\cite{GS1}.
Goldberg's result requires a substitute
for~\eqref{eq:lim_ap} on $L^p(\R^3)$ spaces, rather than weighted $L^2$ spaces. Such a substitute exists, and
is known to be related to the Stein-Tomas theorem in Fourier analysis, see~\cite{stein}. It was first obtained
for the free resolvent by Kenig, Ruiz, and Sogge~\cite{KRS}, and extended to perturbed resolvents by Goldberg
and the author~\cite{GS2}, as well as Ionescu and the author~\cite{IS}. For example, in $\R^3$ the bound from~\cite{KRS}
takes the form
\[ \|R_0(\lambda^2+i\eps)f\|_{L^4(\R^3)} \le C\,\lambda^{-\frac12} \|f\|_{L^{\frac43}(\R^3)}, \]
and in \cite{GS2} it is proved that
\beeq
\label{eq:ST}
\sup_{0<\eps<1,\;\lambda\ge\lambda_0}
\Big\|(-\Laplace+V - (\lambda^2+i\eps))^{-1}\Big\|_{\frac43\to 4} \le C(\lambda_0,V)\;\lambda^{-\half}.
\eneq
for all  real-valued $V\in L^p(\R^3)\cap L^{\frac32}(\R^3), p > \frac32$
and every $\lambda_0>0$. This of course requires absence of imbedded bound states in the continuous spectrum,
which was proved for the same class of $V$ by Ionescu and Jerison~\cite{IonJer}. A very different approach from
the one in~\cite{GS2} to estimates of the form~\eqref{eq:ST} was found in~\cite{IS},
which is related to~\cite{RV}.
\cite{IS} applies to
all dimensions $d\ge2$ and quite general perturbations (including magnetic ones) of $-\Laplace$,
but it also does not rely on~\cite{IonJer}. In fact, as in Agmon's classical paper~\cite{agmon}
 it is shown that the imbedded eigenvalues form a discrete set outside of which a bound as
 in~\eqref{eq:ST} holds (albeit
on somewhat different spaces). Moreover, this is obtained under
the assumption that $V\in L^p(\R^d)$ for some $\frac{d}{2}\le p
\le \frac{d+1}{2}$. The upper limit of $\frac{d+1}{2}$ here is
natural in some ways, since Ionescu and Jerison have found a
smooth, real-valued potential in $L^p(\R^d)$ for all
$p>\frac{d+1}{2}$ which has an imbedded eigenvalue. The lower limit of $d/2$
is the usual one for self-adjointness purposes.

Returning to dispersive estimates, Goldberg~\cite{gold2} proved
that even $V\in L^{p}(\R^3)\cap L^q(\R^3)$ with $p<\frac32<q$
suffices for a dispersive estimate with the usual restriction on
zero energy. Note that this is nearly critical with respect to the
natural scaling of the Schr\"odinger equation in~$\R^3$. One of
his main observations for the low energy argument was that for such $V$ (and assuming that zero is neither an eigenvalue nor a resonance)
\[ \sup_{\lambda\in\R} \big\| (I+VR_0^+(\lambda^2))^{-1} \big\|_{1\to1} < \infty,\]
see~\cite{gold2} for further details.
As far as high energies are concerned, Goldberg noticed that the Born series estimate from~\cite{RS} can be improved so that the
$k^{th}$ term is bounded by $(\lambda_1^{-\eps}\|V\|)^k$ with $\|V\|=\max(\|V\|_p,\|V\|_q)$
as opposed to $(\|V\|_{\kato}/4\pi)^k$. Choosing $\lambda_1$ large this guarantees a
convergent series.

\section{The one-dimensional case}
\label{sec:1d}

We will not repeat the discussion of the one-dimensional theorems from
the introduction where the results from \cite{Wed1}, \cite{Wed2}, \cite{AY}, or~\cite{GS1}
were described. Rather, we would like to focus on a novel estimate that exploits
the absence of a resonance by means of weights and obtains a better rate of decay.
It was motivated by the work of Murata~\cite{Mur} and Buslaev, Perelman~\cite{BP} on improved
local $L^2$ decay in the absence of resonances in dimension one. Note that the weight $\la x\ra$ is optimal in
the sense that it cannot be replaced with $\la x\ra^\tau$, $\tau<1$.

\begin{theorem}
\label{thm:1dimprove}
Suppose $V$ is real-valued and $\|\la x\ra^4 V\|_1<\infty$.
Let $H=-\frac{d^2}{dx^2}+V$ have the property that zero energy is not a resonance. Then
\[ \|\la x\ra^{-1}e^{itH}P_{ac}f \|_\infty \le Ct^{-\frac32} \|\la x\ra f\|_1\]
for all $t>0$.
\end{theorem}
\begin{proof}
Let $\lambda_0 = \|\la x\ra V\|_1^2$ and suppose $\chi$ is a smooth cut-off such that
$\chi(\lambda)=0$ for $\lambda\le \lambda_0$ and $\chi(\lambda)=1$ for $\lambda\ge 2\lambda_0$.
Recall that
\[ R_0(\lambda\pm i0)(x)=\frac{\pm i}{2\sqrt{\lambda}} e^{\pm i|x|\sqrt{\lambda}}. \]
Hence,
\[
\big|\langle R_0(\lambda+i0)(VR_0(\lambda+i0))^n f,g \rangle \big| \le
 (2\sqrt{\lambda})^{-n-1} \|V\|_1^n \|f\|_1\, \|g\|_1, \nn
\]
and the Born series
\begin{equation}
R_V(\lambda\pm i0) = \sum_{n=0}^\infty R_0(\lambda\pm i0)(-VR_0(\lambda\pm i0))^n
\label{eq:born}
\end{equation}
converges  in the operator norm $L^1(\R)\to L^\infty(\R)$ provided $\lambda>\lambda_0$.  The  absolutely
continuous part of the spectral measure is given by
\[
\langle E_{a.c.}(d\lambda) f,g \rangle
= \langle \frac{1}{2\pi i}[R_V(\lambda+i0)-R_V(\lambda-i0)]f,g \rangle \,d\lambda.
\]
Therefore, integrating by parts once yields
\beeq
\label{eq:evol}
\langle e^{itH} \chi(H) f,g \rangle
= -(4\pi t)^{-1}\sum_{n=0}^\infty \int_{-\infty}^\infty e^{it\lambda^2}
\frac{d}{d\lambda}\Big[\chi(\lambda^2) \langle R_0(\lambda^2+i0)
(VR_0(\lambda^2+i0))^n f,g\rangle\Big]\, d\lambda
\eneq
where we have first changed variables $\lambda \to\lambda^2$.
Summation and integration may be exchanged because the Born series
converges absolutely in the $L^1(d\lambda)$ norm, and the domain of integration
is extended to $\R$ via the identity
$R_0(\lambda^2 - i0) = R_0((-\lambda)^2 + i0)$.
The kernel of $R_0(\lambda^2+i0)(VR_0(\lambda^2+i0))^n$ is given explicitly
by the formula
\[ R_0(\lambda^2+i0)(VR_0(\lambda^2+i0))^n(x,y) =
\frac{1}{(2\lambda)^{n+1}}\int_{\R^n} \prod_{j=1}^nV(x_j)
e^{i\lambda(|x-x_1| + |y-x_n| +\sum_{k=2}^n|x_k-x_{k-1}|)} dx \]
with $dx=dx_1\ldots dx_n$.
Hence, in view of the derivative in \eqref{eq:evol},
\begin{align}
&\big|\langle e^{itH} \chi(H) f,g \rangle \big| \nn \\
&\le C |t|^{-1}\sum_{n=0}^\infty (2\sqrt{\lambda_0})^{-n-1}
\sup_{a\in\R} \left| \int_{-\infty}^\infty e^{i(t\lambda^2+a\lambda)}
\chi(\lambda^2)\,\lambda^{-n-1}\lambda_0^{(n+1)/2}\, d\lambda \right| \|\la x\ra V\|_1^n \,
\|\la x\ra f\|_1 \|\la x\ra g\|_1
\nn \\
&\quad + C |t|^{-1}\sum_{n=0}^\infty (2\sqrt{\lambda_0})^{-n-1}
\sup_{a\in\R} \left| \int_{-\infty}^\infty e^{i(t\lambda^2+a\lambda)}
\chi'(\lambda^2)\,\lambda^{-n}\lambda_0^{(n+1)/2}\, d\lambda \right| \|V\|_1^n \, \|f\|_1 \|g\|_1
\nn \\
&\le C(V)\, |t|^{-\frac32} \|\la x\ra f\|_1\|\la x\ra g\|_1. \label{eq:decay2}
\end{align}
We used the dispersive bound for the one-dimensional Schr\"odinger
equation to pass to~\eqref{eq:decay2}, observing in particular that
\[ \sup_{n\ge0} \big\| [\chi(\lambda^2)\,\lambda^{-n-1}\lambda_0^{(n+1)/2}]^{\vee} \big\| < \infty \]
where the norm refers to the total variation norm of measures.

It remains to consider small energies, i.e., those $\lambda$ for which $\chi(\lambda^2)\ne1$.
In this case, we let $f_j(\cdot,\lambda)$ for $j=1,2$ be the Jost solutions. They satisfy
\[ \big(-\frac{d^2}{dx^2}+V-\lambda^2\big)f_j(x,\lambda)=0,\; f_1(x,\lambda)\sim e^{ix\lambda}
\text{\ \ as\ \ }x\to\infty,\;
f_2(x,\lambda)\sim e^{-ix\lambda} \text{\ \ as\ \ }x\to-\infty
\]
for any $\lambda\in\R$. Furthermore, if $\lambda\ne0$, then
\beeq
\label{eq:RT}
 f_1(\cdot,\lambda)=\frac{R_1(\lambda)}{T(\lambda)}f_2(\cdot,\lambda)+\frac{1}{T(\lambda)}f_2(\cdot,-\lambda),
\quad f_2(\cdot,\lambda)=
\frac{R_2(\lambda)}{T(\lambda)}f_1(\cdot,\lambda)+\frac{1}{T(\lambda)}f_1(\cdot,-\lambda),
\eneq
where $T(\lambda)=\frac{-2i\lambda}{W(\lambda)}$ with $W(\lambda)=W[f_1(\cdot,\lambda),f_2(\cdot,\lambda)]$ and
\[ R_1(\lambda)=-\frac{T(\lambda)}{2i\lambda}W[f_1(\cdot,\lambda),f_2(\cdot,-\lambda)], \;
 R_2(\lambda)=\frac{T(\lambda)}{2i\lambda}W[f_1(\cdot,-\lambda),f_2(\cdot,\lambda)].
\]
Then the jump condition of the resolvent $R_V$ across the spectrum takes the form
\[ \big(R_V(\lambda^2+i0)-R_V(\lambda^2-i0)\big)(x,y) =
\frac{|T(\lambda)|^2}{-2i\lambda}(f_1(x,\lambda)f_1(y,-\lambda)+
f_2(x,\lambda)f_2(y,-\lambda))
\]
with $\lambda\ge0$.
Let us denote the distorted Fourier basis by
\[ e(x,\lambda)= \frac{1}{\sqrt{2\pi}}
\left\{ \begin{array}{ll}
T(\lambda) f_1(\cdot,\lambda) & \text{\ \ if\ }\lambda \ge0 \\
T(-\lambda) f_2(x,-\lambda) &  \text{\ \ if\ }\lambda < 0
\end{array} \right.
\]
see Weder's papers~\cite{Wed1} and~\cite{Wed2} for more details on this basis.
Then the evolution $e^{itH}(1-\chi(H))P_{a.c.}$ can be written as
\begin{align}
& \big\la e^{itH}(1-\chi(H))P_{a.c.}\phi,\psi\big\ra \nn \\
&= \frac{1}{2\pi i}\int_{\R^2}\int_0^\infty 2\lambda e^{it\lambda^2}(1-\chi(\lambda^2))
\big(R_V(\lambda^2+i0)-R_V(\lambda^2-i0)\big)(x,y)\,d\lambda \,
\bar{\phi}(x)\psi(y)\, dxdy \nn \\
&= \int_{-\infty}^\infty e^{it\lambda^2} (1-\chi(\lambda^2))\la \psi,e(\cdot,\lambda)\ra
\la e(\cdot,\lambda),\phi\ra \, d\lambda. \label{eq:e+rep}
\end{align}
Our assumption that zero energy is not a resonance implies that $T(\lambda)=\alpha\lambda+o(\lambda)$
where $\alpha\ne0$.
In particular, $T(0)=0$ and $R_1(0)=R_2(0)=-1$.
Integrating by parts in~\eqref{eq:e+rep} therefore yields
\begin{align}
& \big\la e^{itH}(1-\chi(H))P_{a.c.}\phi,\psi\big\ra  \nn \\
& = - \frac{1}{4\pi it}\int_0^\infty e^{it\lambda^2}
\partial_\lambda\Big[(1-\chi(\lambda^2))|T(\lambda)|^2 \lambda^{-1}
\la \psi,f_1(\cdot,\lambda)\ra
\la f_1(\cdot,\lambda),\phi\ra \Big] \, d\lambda \label{eq:f1teil}\\
& \quad - \frac{1}{4\pi it}\int_{-\infty}^0 e^{it\lambda^2}
\partial_\lambda\Big[(1-\chi(\lambda^2))|T(\lambda)|^2 \lambda^{-1}
\la \psi,f_2(\cdot,-\lambda)\ra
\la f_2(\cdot,-\lambda),\phi\ra \Big]\, d\lambda.\nn
\end{align}
By symmetry, it will suffice to treat the integral involving $f_1(\cdot,\lambda)$.
We distinguish three cases, depending on where the derivative $\partial_\lambda$ falls.
We start with the integral
\beeq
\label{eq:part1}
 \int_0^\infty e^{it\lambda^2} \omega(\lambda) f_1(x,\lambda)f_1(y,-\lambda)\, d\lambda,
\eneq
where we have set $\omega(\lambda)=\partial_\lambda[(1-\chi(\lambda^2))|T(\lambda)|^2 \lambda^{-1}]$.
By the preceding, $\omega$ is a smooth function with compact support in $[0,\infty)$.
As usual, we will estimate \eqref{eq:part1} by means of a Fourier transform in $\lambda$. Since we
are working on a half-line, this will actually be a cosine transform. Let $\tilde\omega$ be another
cut-off function satisfying $\omega\tilde\omega=\omega$. Then
\beeq
 \Big| \int_0^\infty e^{it\lambda^2} \omega(\lambda) f_1(x,\lambda)f_1(y,-\lambda)\, d\lambda \Big|
 \le C|t|^{-\half} \| \,[\omega f_1(x,\cdot)]^{\vee} \|_1 \|\, [\tilde\omega f_1(y,-\cdot)]^{\vee} \|_1.
\label{eq:2L1}
\eneq
It remains to estimate
\beeq
\label{eq:omf1}
[\omega f_1(x,\cdot)]^{\vee}(u) := \int_0^\infty \cos(u\lambda) \omega(\lambda) f_1(x,\lambda)\, d\lambda
\eneq
in $L^1$ relative to $u$. The second $L^1$-norm in \eqref{eq:2L1} is treated the same way.
We need to consider the cases $x\ge0$ and $x\le0$ separately. In the former case,
\begin{align}
& [\omega f_1(x,\cdot)]^{\vee}(u) := \int_0^\infty \cos(u\lambda)
e^{ix\lambda} \omega(\lambda) e^{-ix\lambda}f_1(x,\lambda)\, d\lambda \nn \\
&= \half\int_0^\infty e^{i(x+u)\lambda} \omega(\lambda) e^{-ix\lambda}f_1(x,\lambda)\, d\lambda +
\half\int_0^\infty e^{i(x-u)\lambda} \omega(\lambda) e^{-ix\lambda}f_1(x,\lambda)\, d\lambda.
\label{eq:upmx}
\end{align}
If $||u|-|x||\le |x|$, then we simply estimate
\[ |[\omega f_1(x,\cdot)]^{\vee}(u)| \le C.\]
On the other hand, if $||u|-|x||>|x|$, then we integrate by parts in \eqref{eq:upmx}:
\begin{align}
 [\omega f_1(x,\cdot)]^{\vee}(u) &= -\frac{1}{2i(x+u)} \omega(0)f_1(x,0) -
 \frac{1}{2i(x-u)} \omega(0)f_1(x,0) \label{eq:byparts1}\\
&  -\frac{1}{2i(x+u)}\int_0^\infty e^{i(x+u)\lambda} \partial_\lambda\Big[\omega(\lambda) e^{-ix\lambda}f_1(x,\lambda)\Big]\, d\lambda \nn\\
& - \frac{1}{2i(x-u)} \int_0^\infty e^{i(x-u)\lambda} \partial_\lambda\Big[\omega(\lambda) e^{-ix\lambda}f_1(x,\lambda)\Big]\, d\lambda. \nn
\end{align}
Since
\[ \sup_{x\ge0,\,\lambda}| \partial^j_\lambda[\omega(\lambda) e^{-ix\lambda}f_1(x,\lambda)]|\le C(V),\]
for $j=0,1,2$,  it follows that
\[
|[\omega f_1(x,\cdot)]^{\vee}(u)| \le C \frac{|x|}{|x^2-u^2|} + C(u+x)^{-2} + C(u-x)^{-2}.
\]
The conclusion is that
\beeq
\label{eq:xest1} \int_{\R} |[\omega f_1(x,\cdot)]^{\vee}(u)|\, du \le C\la x\ra.
\eneq
To deal with $x\le0$, we use \eqref{eq:RT}. Thus,
\begin{align}
 [\omega f_1(x,\cdot)]^{\vee}(u) &= \int_0^\infty \cos(u\lambda) \omega(\lambda)\frac{R_1(\lambda)+1}{T(\lambda)} f_2(x,\lambda)\, d\lambda
\label{eq:Rint}\\
& \quad + \int_0^\infty \cos(u\lambda)\omega(\lambda)\frac{1}{T(\lambda)} (f_2(x,\lambda)-f_2(x,-\lambda))\, d\lambda. \label{eq:Tint}
\end{align}
Set $\omega_1=\omega(\lambda)\frac{R_1(\lambda)+1}{T(\lambda)}$. Then \eqref{eq:Rint} can be written as
\[ \int_0^\infty \cos(u\lambda) \omega(\lambda)\frac{R_1(\lambda)+1}{T(\lambda)} f_2(x,\lambda)\, d\lambda =
\int_0^\infty \cos(u\lambda)e^{-ix\lambda} \omega_1(\lambda) e^{ix\lambda} f_2(x,\lambda)\, d\lambda.
\]
Hence it can be treated by the same arguments as \eqref{eq:omf1} with $x\ge0$. Indeed, simply use that
\[ \sup_{x\le0,\,\lambda}| \partial_\lambda[\omega_1(\lambda) e^{ix\lambda}f_2(x,\lambda)]|\le C(V). \]
On the other hand, \eqref{eq:Tint} is the same as (with $\partial_2$ being the partial derivative with
respect to the second variable of $f_2$)
\begin{align}
&\int_{-1}^1 \int_0^\infty \cos(u\lambda)\omega(\lambda)\frac{\lambda}{T(\lambda)} \partial_2 f_2(x,\lambda\sigma)\, d\lambda d\sigma \nn \\
&= \int_{-1}^1 \int_0^\infty \cos(u\lambda)e^{-i\lambda x\sigma}\omega_2(\lambda) \partial_2[ e^{i\lambda x\sigma} f_2(x,\lambda\sigma)]\, d\lambda d\sigma \label{eq:1inte}\\
& \quad - ix\int_{-1}^1 \int_0^\infty \cos(u\lambda)e^{-i\lambda x\sigma}\omega_2(\lambda)  e^{i\lambda x\sigma} f_2(x,\lambda\sigma)\, d\lambda d\sigma,  \label{eq:2inte}
\end{align}
where we have set $\omega_2(\lambda)=\omega(\lambda)\frac{\lambda}{T(\lambda)}$ (a smooth, compactly supported function in $[0,\infty)$).
We will focus on the second integral \eqref{eq:2inte}, since the first one~\eqref{eq:1inte} is similar.
We will integrate by parts in $\lambda$, but only on the set $|\sigma x\pm u|\ge1$. Then
\begin{align}
&- ix\int_{-1}^1 \int_0^\infty \cos(u\lambda)e^{-i\lambda x\sigma}\omega_2(\lambda)  e^{i\lambda x\sigma} f_2(x,\lambda\sigma)\, d\lambda
 \chi_{[|\sigma x\pm u|\ge1]}\,   d\sigma \nn \\
&= \int_{-1}^1\frac{x}{2(-\sigma x+u)} \omega_2(0)f_2(x,0)\, \chi_{[|\sigma x\pm u|\ge1]}\, d\sigma + \int_{-1}^1\frac{x}{2(-\sigma x-u)} \omega_2(0)f_2(x,0) \, \chi_{[|\sigma x\pm u|\ge1]}\, d\sigma \nn \\
&  +\int_{-1}^1\frac{x}{2(-\sigma x+u)}\int_0^\infty e^{i(-\sigma x+u)\lambda} \partial_\lambda\Big[\omega_2(\lambda) e^{ix\lambda\sigma}f_2(x,\lambda)\Big]\, d\lambda \, \chi_{[|\sigma x\pm u|\ge1]}\, d\sigma \label{eq:gut1}\\
& + \int_{-1}^1 \frac{x}{2(-\sigma x-u)} \int_0^\infty e^{i(-\sigma x-u)\lambda} \partial_\lambda\Big[\omega_2(\lambda) e^{ix\lambda\sigma}f_2(x,\lambda)\Big]\, d\lambda\,  \chi_{[|\sigma x\pm u|\ge1]}\,  d\sigma. \label{eq:gut2}
\end{align}
The first two integrals here (which are due to the boundary $\lambda=0$) contribute
\[
\int_{-1}^1\frac{x}{2(-\sigma x+u)} \omega_2(0)f_2(x,0)\, \chi_{[|\sigma x\pm u|\ge1]}\, d\sigma + \int_{-1}^1\frac{x}{2(\sigma x-u)} \omega_2(0)f_2(x,0) \, \chi_{[|\sigma x\pm u|\ge1]}\, d\sigma =0,
\]
where we performed a change of variables $\sigma\mapsto -\sigma$ in the second one.
Integrating by parts one more time in \eqref{eq:gut1} and \eqref{eq:gut2} with respect to $\lambda$ implies
\begin{align*}
& \int_{-\infty}^\infty \Big|x\int_{-1}^1 \int_0^\infty \cos(u\lambda)e^{-i\lambda x\sigma}\omega_2(\lambda)  e^{i\lambda x\sigma} f_2(x,\lambda\sigma)\, d\lambda
 \chi_{[|\sigma x\pm u|\ge1]}\,   d\sigma\Big|\, du \nn \\
& \le C \int_{-\infty}^\infty \int_{-1}^1\frac{|x|}{(-\sigma x+u)^2} \chi_{[|\sigma x\pm u|\ge1]}\, d\sigma du
+ C \int_{-\infty}^\infty \int_{-1}^1\frac{|x|}{(-\sigma x-u)^2} \chi_{[|\sigma x\pm u|\ge1]}\, d\sigma du \le C\, |x|.
\end{align*}
Finally, the cases $|\sigma x+u|\le1 $ and $|\sigma x-u|\le1$ each contribute at most $C|x|$ to the $u$-integral. Hence
\[
 \int_{-\infty}^\infty \Big|x\int_{-1}^1 \int_0^\infty \cos(u\lambda)e^{-i\lambda x\sigma}\omega_2(\lambda)  e^{i\lambda x\sigma} f_2(x,\lambda\sigma)\, d\lambda d\sigma \Big|\, du
 \le C|x|.
\]
Since \eqref{eq:1inte} can be treated the same way (in fact, the bound is $O(1)$), we obtain
\[
\int_{-\infty}^\infty\Big| \int_{-1}^1 \int_0^\infty \cos(u\lambda)\omega(\lambda)\frac{\lambda}{T(\lambda)} \partial_2 f_2(x,\lambda\sigma)\, d\lambda d\sigma
\Big| \, du \le C\la x\ra.
\]
In view of \eqref{eq:xest1}, \eqref{eq:Rint}, and \eqref{eq:Tint},
\[
 \big\| [\omega f_1(x,\cdot)]^{\vee} \big\|_1 \le C\la x\ra \quad \forall \;x\in\R,
\]
which in turn implies that
\beeq
\label{eq:uff1}
\Big| \int_0^\infty e^{it\lambda^2} \omega(\lambda) f_1(x,\lambda)f_1(y,-\lambda)\, d\lambda \Big| \le C\,|t|^{-\half} \la x\ra\la y\ra,
\eneq
for all $x,y\in \R$, see \eqref{eq:part1}.
This is the desired estimate on \eqref{eq:f1teil}, but only for the case when $\partial_\lambda$ falls on the factors not involving $f_1$.
We now consider the case when $\partial_\lambda$ falls on $f_1(x,\lambda)$. The integral in which $\partial_\lambda$ falls on
$f_1(y,-\lambda)$ is analogous. Hence, we need to estimate
\begin{align}
&\int_0^\infty e^{it\lambda^2} (1-\chi(\lambda^2))|T(\lambda)|^2 \lambda^{-1}
\partial_\lambda f_1(x,\lambda)
\;f_1(y,-\lambda) \, d\lambda \nn \\
& = \int_0^\infty e^{it\lambda^2} (1-\chi(\lambda^2))T(-\lambda) \lambda^{-1}
\partial_\lambda [T(\lambda)f_1(x,\lambda)]
\, f_1(y,-\lambda) \, d\lambda \label{eq:neu} \\
& + \int_0^\infty e^{it\lambda^2} (1-\chi(\lambda^2))T(-\lambda)T'(\lambda)\lambda^{-1}  f_1(x,\lambda)
\;f_1(y,-\lambda) \, d\lambda. \label{eq:zuruck}
\end{align}
The integral in \eqref{eq:zuruck} is of the same form as that in \eqref{eq:part1}.
It therefore suffices to control~\eqref{eq:neu}.  Let $\omega_3(\lambda)=(1-\chi(\lambda^2))T(-\lambda)\lambda^{-1}$.
By the same reductions as before, we need to show that
\[
 \big\| \big[\omega_3 \partial_\lambda[T(\lambda)f_1(x,\cdot)]\big]^{\vee} \big\|_1 \le C\la x\ra \quad \forall \;x\in\R.
\]
Thus consider
\begin{align*}
 \int_0^\infty \cos(u\lambda)\omega_3(\lambda) \partial_\lambda[T(\lambda)f_1(x,\lambda)]\, d\lambda
& = ix\int_0^\infty \cos(u\lambda) e^{ix\lambda}\omega_3(\lambda) T(\lambda)e^{-ix\lambda}f_1(x,\lambda)\, d\lambda  \\
& + \int_0^\infty \cos(u\lambda) e^{ix\lambda}\omega_3(\lambda) \partial_\lambda[ T(\lambda)e^{-ix\lambda}f_1(x,\lambda)]\, d\lambda.
\end{align*}
If $x\ge0$, integrating by parts leads to
\begin{align}
& ix\int_0^\infty \cos(u\lambda) e^{ix\lambda}\omega_3(\lambda) T(\lambda)e^{-ix\lambda}f_1(x,\lambda)\, d\lambda   \label{eq:odin}\\
& = -\frac{ix}{2i(x+u)}\int_0^\infty e^{i(x+u)\lambda} \partial_\lambda\Big[\omega_3(\lambda)T(\lambda) e^{-ix\lambda}f_1(x,\lambda)\Big]\, d\lambda \nn\\
&\quad - \frac{ix}{2i(x-u)} \int_0^\infty e^{i(x-u)\lambda} \partial_\lambda\Big[\omega_3(\lambda)T(\lambda) e^{-ix\lambda}f_1(x,\lambda)\Big]\, d\lambda\nn
\end{align}
as well as
\begin{align}
& \int_0^\infty \cos(u\lambda) e^{ix\lambda}\omega_3(\lambda) \partial_\lambda[ T(\lambda)e^{-ix\lambda}f_1(x,\lambda)]\, d\lambda \label{eq:dwa} \\
&= -\frac{1}{2i(x+u)} \omega_3(0)\partial_\lambda[ T(\lambda)e^{-ix\lambda}f_1(x,\lambda)]\Big\vert_{\lambda=0} -
\frac{1}{2i(x-u)} \omega_3(0) \partial_\lambda[ T(\lambda)e^{-ix\lambda}f_1(x,\lambda)]\Big\vert_{\lambda=0} \nn\\
& \quad -\frac{1}{2i(x+u)}\int_0^\infty e^{i(x+u)\lambda} \partial_\lambda\Big[\omega_3(\lambda) \partial_\lambda[ T(\lambda)e^{-ix\lambda}f_1(x,\lambda)\Big]\, d\lambda \nn \\
& \quad- \frac{1}{2i(x-u)} \int_0^\infty e^{i(x-u)\lambda} \partial_\lambda\Big[\omega_3(\lambda) \partial_\lambda[ T(\lambda)e^{-ix\lambda}f_1(x,\lambda)\Big]\, d\lambda. \nn
\end{align}
Integrating by parts one more time in \eqref{eq:odin} implies
\begin{align*}
& \Big|ix\int_0^\infty \cos(u\lambda) e^{ix\lambda}\omega_3(\lambda) T(\lambda)e^{-ix\lambda}f_1(x,\lambda)\, d\lambda \Big|  \nn \\
& \le C|x|(1+|x-u|)^{-2}+C|x|(1+|x+u|)^{-2}
\end{align*}
uniformly in $x\ge0$, whereas \eqref{eq:dwa} is treated the same way as~\eqref{eq:byparts1}.
One needs to use here that
\[ \sup_{x\ge0,\,\lambda}| \partial^j_\lambda[\omega_3(\lambda) e^{-ix\lambda}f_1(x,\lambda)]|\le C(V),\]
for $j=0,1,2,3$ which follows from $\|\la x\ra^4 V\|_1<\infty$.
Consequently, we have proved that
\[ \int_{\R}\Big|\int_0^\infty \cos(u\lambda)\omega_3(\lambda) \partial_\lambda[T(\lambda)f_1(x,\lambda)]\, d\lambda \Big|\,du\le C\la x\ra\]
uniformly in $x\ge0$.

Next, we deal with the case $x\le0$.
In view of~\eqref{eq:RT},
\[ T(\lambda) f_1(\cdot,\lambda) = R_1(\lambda) f_2(\cdot,\lambda) + f_2(\cdot,-\lambda).\]
This implies that
\begin{align}
& \int_0^\infty \cos(u\lambda)\omega_3(\lambda) \partial_\lambda[T(\lambda)f_1(x,\lambda)]\, d\lambda \label{eq:dTf1}\\
& = \int_0^\infty \cos(u\lambda)e^{-ix\lambda}\omega_3(\lambda) \partial_\lambda[R_1(\lambda)e^{ix\lambda}f_2(x,\lambda)]\, d\lambda
 -ix\int_0^\infty \cos(u\lambda)e^{-ix\lambda} \omega_3(\lambda) R_1(\lambda)e^{ix\lambda}f_2(x,\lambda)\, d\lambda  \nn \\
& + \int_0^\infty \cos(u\lambda)e^{-ix\lambda}\omega_3(\lambda) \partial_\lambda[e^{ix\lambda}f_2(x,-\lambda)]\, d\lambda
+ ix\int_0^\infty \cos(u\lambda)e^{ix\lambda} \omega_3(\lambda) e^{-ix\lambda}f_2(x,-\lambda)\, d\lambda. \nn
\end{align}
The two integrals which are not preceded by factors of $ix$ are treated just as in~\eqref{eq:dwa}. The only difference
here is that the estimates are uniform in $x\le0$ rather than $x\ge0$.
On the other hand, the integrals preceded by $ix$ need to be  integrated by parts in $\lambda$. It is important to check
that the boundary terms at $\lambda=0$ do not contribute to this case. Indeed, these boundary terms are
\begin{align*}
& \frac{x}{2(u-x)}\omega_3(0)R_1(0)f_2(x,0)-\frac{x}{2(u+x)}\omega_3(0)R_1(0)f_2(x,0) \\
& -\frac{x}{2(u+x)}\omega_3(0)f_2(x,0)
-\frac{x}{2(x-u)}\omega_3(0)f_2(x,0)=0,
\end{align*}
since $R_1(0)=-1$. Hence, integrating by parts leads to an expression similar to~\eqref{eq:odin}. The conclusion is that \eqref{eq:dTf1}
satisfies
\[\int_{\R} \Big|\int_0^\infty \cos(u\lambda)\omega_3(\lambda) \partial_\lambda[T(\lambda)f_1(x,\lambda)]\, d\lambda \Big|\,du\le C\la x\ra\]
uniformly in $x\le 0$, and we are done.
\end{proof}

In~\cite{KS} the same bound is proved for non-selfadjoint systems of the type that arise
by linearizing NLS around a ground state standing wave. It is crucial for proving the existence
of stable manifolds for all super-critical NLS in one dimension.

In dimension one, there is some recent work of Cai~\cite{cai2} on
dispersion for Hill's operator. More precisely, let
$H=-\frac{d^2}{dx^2} + q$ where $q$ is periodic and such that its
spectrum has precisely one gap. It is well-known that such $q$ are
characterized in terms of Weierstrass elliptic functions. As part
of his Caltech Ph.D.~thesis,  Cai showed that for this $H$ one
always has
\[ \|e^{itH}f\|_\infty \le Ct^{-\frac14}\|f\|_1, \qquad t\ge1\]
and that generically in the potential one can replace $\frac14$ with $\frac13$.

\section{The two-dimensional case}

The following two-dimensional dispersive estimate was obtained in~\cite{Sch}.

\begin{theorem}
\label{thm:main2d}
Let $V:\R^2\to\R$ be a measurable function
such that $|V(x)|\le C(1+|x|)^{-\beta}$, $\beta>3$. Assume in addition
that zero is a regular point of the spectrum of $H=-\Laplace+V$. Then
\[ \big\|e^{itH}P_{ac}(H) f\big \|_{\infty} \le C|t|^{-1}\|f\|_1\]
for all $f\in L^1(\R^2)$.
\end{theorem}

The definition of zero being a regular point amounts to the following, see
Jensen, Nenciu~\cite{JenNen}:
{\em Let $V\not\equiv 0$ and set $U={\mathrm sign}\, V$, $v=|V|^{\half}$.
Let $P_v$ be the orthogonal projection onto $v$ and set $Q=I-P_v$. Finally, let
\[ (G_0 f)(x):= -\frac{1}{2\pi}\int_{\R^2} \log|x-y|\, f(y)\, dy.\]
Then zero is regular iff $Q(U+vG_0v)Q$ is invertible on $QL^2(\R^2)$.}

Jensen and Nenciu study $\ker[Q(U+vG_0v)Q]$ on $QL^2(\R^2)$.
It can be completely described in terms of
solutions~$\Psi$ of $H\Psi=0$. In particular, its dimension is at most three plus the dimension of the
zero energy eigenspace, see Theorem~6.2 and Lemma~6.4  in~\cite{JenNen}.  The extra three
dimensions here are called resonances.
Hence, the requirement that zero is a regular point is the analogue of the usual
condition that zero is neither an eigenvalue nor a resonance of~$H$. An equivalent
characterization of a regular point was given in~\cite{BGD2},
albeit under the additional assumption that $\int_{\R^2} V(x)\, dx\ne0$.

As far as the spectral properties of $H$ are concerned, we note that under the hypotheses
of Theorem~\ref{thm:main2d} the spectrum of $H$ on $[0,\infty)$ is purely absolutely
continuous, and that the spectrum is pure point on $(-\infty,0)$ with at most finitely
many eigenvalues of finite multiplicities. The latter follows for example from Stoiciu~\cite{mihai},
who obtained Birman-Schwinger type bounds in the case of two dimensions.

Theorem~\ref{thm:main2d} appears to be the first $L^1\to L^\infty$ bound
with $|t|^{-1}$ decay in~$\R^2$. Yajima~\cite{Y4} and Jensen, Yajima~\cite{JenYaj}
proved the $L^p(\R^2)$ boundedness of the wave operators under stronger decay
assumptions on~$V(x)$, but only for $1<p<\infty$. Hence their result does
not imply Theorem~\ref{thm:main2d}. Local $L^2$ decay was studied by
Murata~\cite{Mur}, but he does not consider $L^1\to L^\infty$ estimates.

The main challenge in two dimensions is of course the low energy
part. This is due to the fact that the free resolvent
$R_0^{\pm}(\lambda^2)=(-\Laplace-(\lambda^2\pm i0))^{-1}$ has the
kernel ($H_0^{\pm}$ being the Hankel functions)
\[  R_0^{\pm}(\lambda^2)(x,y)= \pm \frac{i}{4} H_0^{\pm}(\lambda|x-y|), \]
which is singular at energy zero (which, just
as in dimension one, expresses the fact that the free problem has a resonance at zero).
It is a consequence of the asymptotic expansion of Hankel functions that for all $\lambda>0$,
\begin{equation}
\label{eq:zero_exp}
R_0^{\pm}(\lambda^2)  = \Big[\pm\frac{i}{4} - \frac{1}{2\pi}\gamma - \frac{1}{2\pi} \log(\lambda/2)\Big]P_0
+ G_0  + E_0^{\pm}(\lambda).
\end{equation}
Here $P_0f:=\int_{\R^2} f(x)\,dx$,
$G_0f(x)=-\frac{1}{2\pi}\int_{\R^2} \log|x-y|\,f(y)\,dy$, and the
error $E_0^{\pm}(\lambda)$ has the property that \beeq
\label{eq:err}
\big\|\sup_{0<\lambda}\lambda^{-\half}|E_0^{\pm}(\lambda)(\cdot,\cdot)|\,\big
\|+ \big\|\sup_{0<\lambda}\lambda^{\half} |\partial_\lambda
E_0^{\pm}(\lambda)(\cdot,\cdot)|\,\big\| \les 1 \eneq with respect
to the Hilbert-Schmidt norm in
$\bdd(L^{2,s}(\R^2),L^{2,-s}(\R^2))$ with $s>\frac32$. These error
estimates may seem artificial, but they allow for the least amount
of decay on~$V$. The following lemma from~\cite{Sch} contains the
expansion of the perturbed resolvent around energy zero needed in
the proof of Theorem~\ref{thm:main2d}. It displays an important
idea from~\cite{JenNen}, namely to re-sum infinite series of
powers of $\log\lambda$ into one function $h_{\pm}(\lambda)$. This
feature is crucial for our purposes. Given $V\not\equiv 0$, set
$U={\rm sign}\, V$, $v=|V|^{\half}$. Let $P_v$ be the orthogonal
projection onto $v$ and set $Q=I-P_v$. Finally, let
$D_0=[Q(U+vG_0v)Q]^{-1}$ on $QL^2(\R^2)$

\begin{lemma}
\label{lem:inv_exp}
Suppose that zero is a regular point of the spectrum of $H=-\Laplace+V$. Then
for some sufficiently small $\lambda_1>0$, the operators
$M^{\pm}(\lambda) :=U+vR_0^{\pm}(\lambda^2)v $ are invertible
for all $0<\lambda<\lambda_1$ as bounded operators on $L^2(\R^2)$, and one has the expansion
\beeq
\label{eq:M_inv_exp}
M^{\pm}(\lambda)^{-1} = h_{\pm}(\lambda)^{-1}S+QD_0Q+E^{\pm}(\lambda),
\eneq
where $h_{+}(\lambda)=a\log \lambda+z$, $a$ is real, $z$ complex, $a\ne0$, $\Im z\ne0$, and $h_{-}(\lambda)=\overline{h_+(\lambda)}$.
Moreover, $S$ is of finite rank and has a real-valued kernel,
and $E^{\pm}(\lambda)$ is a Hilbert-Schmidt operator that satisfies the bound
\beeq
\label{eq:Eest}
\big \|\sup_{0<\lambda<\lambda_1} \lambda^{-\half} |E^{\pm}(\lambda)(\cdot,\cdot)|\,\big  \|_{HS}
+  \big \|\sup_{0<\lambda<\lambda_1} \lambda^{\half} |\partial_\lambda E^{\pm}(\lambda)(\cdot,\cdot)|\, \big \|_{HS} \les 1
\eneq
where the norm refers to the Hilbert-Schmidt norm on $L^2(\R^2)$. Finally, let $R_V^{\pm}(\lambda^2) = (-\Laplace+V-(\lambda^2\pm i0))^{-1}$.  Then
\beeq
R_V^{\pm}(\lambda^2)  = R_0^{\pm}(\lambda^2) -
R_0^{\pm}(\lambda^2)v M^{\pm}(\lambda)^{-1} v R_0^{\pm}(\lambda^2).
\label{eq:RV}
\eneq
This is to be understood as an identity between operators $L^{2,\half+\eps}(\R^2)\to L^{2,-\half-\eps}(\R^2)$ for some sufficiently small $\eps>0$.
\end{lemma}

The low energy part of the proof of Theorem~\ref{thm:main2d} is based on a careful estimation of the contribution of
each of the terms in~\eqref{eq:M_inv_exp} to $R_V$ in~\eqref{eq:RV} by means of the method of stationary phase, see~\cite{Sch}.

Murata~\cite{Mur} discovered that under the assumptions of Theorem~\ref{thm:main2d}
\[
\big\| w e^{itH}P_{ac}(H) w\,f\big \|_{2} \le C|t|^{-1}(\log t)^{-2}\|f\|_2
\]
provided $w(x)=\la x\ra^{-\sigma}$ with some sufficiently large
$\sigma>0$. In other words, he obtained improved local $L^2$ decay
provided zero energy is regular. Needless to say, such improved
decay is impossible for the  $L^1\to L^\infty$ bound, but a
weighted $L^1\to L^\infty$ estimate as in
Theorem~\ref{thm:1dimprove} with the improved $|t|^{-1}(\log
t)^{-2}$ decay is quite possibly true but currently unknown. Due
to the integrability of this decay at infinity, such a bound would
be useful for the study of nonlinear asymptotic stability of
(multi) solitons in dimension two.

\section{Time-dependent potentials}

It seems unreasonable to expect a general theory of dispersion for the Schr\"odinger equation
\beeq
\label{eq:Vt}
 i\partial_t \psi+ \Laplace \psi + V(t,\cdot)\psi =0
\eneq
for time-dependent potentials $V(t,\cdot)$. While the $L^2$ norm is preserved for real-valued $V$, it is
well-known that in contrast to time-independent $V$ higher $H^s$ norms can grow in this case, see e.g.\ Bourgain~\cite{master1}, \cite{master2},
and Erdogan, Killip, Schlag~\cite{EKS}.

The classical work of Davies~\cite{Da}, Howland~\cite{Howl1}, \cite{Howl2}, \cite{Howl3},  and Yajima~\cite{Y0},
deals with scattering and wave operators in this context.
Recall that if $U(t,s)$ denotes the evolution of~\eqref{eq:Vt} from time $s$ to time~$t$, then
\[ W_{\pm}(s) = s-\lim_{t\to\pm\infty} e^{-i(t-s)\Laplace}U(t,s) \]
are the wave operators (strictly speaking, the existence of these limits is usually refered to as {\em completeness}, but
we are following Howland's terminology).
In  analogy to the treatment of time-dependent Hamiltonians in classical mechanics,
Howland~\cite{Howl2} develops a formalism for treating
time-dependent potentials in which $K=-i\partial_t+H(t)$ is considered as a self-adjoint operator on the Hilbert space $L^2(-\infty,\infty;L^2(\R^d))$.
He shows that the existence of $W_{\pm}$ is equivalent to the existence of the strong limits
\[ \calW_{\pm}:=s-\lim_{\sigma\to\pm\infty} e^{i\sigma \calK_0}e^{-i\sigma \calK} \]
and that $\calW_{\pm}$ is the same as multiplication by $W_{\pm}(t)$. Furthermore, following Kato~\cite{kato},
he formulates a condition which insures that the wave operators are unitary. He applies this to~\eqref{eq:Vt}
with (real-valued) potentials
\[ V \in L^{r+\eps}_t(L^p_x)\cap L^{r-\eps}_t(L^p_x), \quad r=\frac{2p}{2p-d}, \; \frac{d}{2}<p\le \infty,\; d>1\]
to conclude that for such $V$ the wave operators exist and are unitary.
In~\cite{Howl2}, Howland obtained similar results for $d\ge3$ potentials that are small at infinity (rather than
vanishing).

Dispersive estimates were obtained by Rodnianksi and the author~\cite{RS} for small but not necessarily
decaying time-dependent potentials in $\R^3$, whereas the case of decaying $V$ and dimensions $\ge2$
was studied by Naibo, Stepanov~\cite{NS}, and
d'Ancona, Pierfelice, Visciglia~\cite{APV}. In particular, the result from~\cite{RS} insures that in $\R^3$ and for small $\eps$
\[  i\partial_t \psi+ \Laplace \psi + \eps F(t)V(x)\psi = 0 \]
has the usual $t^{-\frac32}$ dispersive $L^1\to L^\infty$ decay for any real-valued trigonometric polynomial $F(t)$ (or more generally,
any quasi-periodic analytic function $F(t)$) and $V$ satisfying $\|V\|_{\kato}<\infty$, see~\eqref{eq:kato}.

Another much studied case is that of {\em time-periodic} $V$, see \cite{Da}, \cite{Howl3}, and~\cite{Y0}.
Suppose $T>0$ is the smallest period of $V$. Then the theory of~\eqref{eq:Vt} reduces to that of the Floquet operator $\calU=U(T,0)$.
The Floquet operator can exhibit bound states and the question arises as to the existence and ranges of the wave operators
(the so called completeness problem) as well as the structure of the discrete spectrum. These issues are addressed in
the aforementioned  references.

More recently, in~\cite{GJY}, Galtbayar, Jensen, and Yajima show that on the orthogonal complement of
the bound states of the Floquet operator the solutions decay locally in~$L^2(\R^3)$.
In addition, O.~Costin, R.~Costin, Lebowitz, and Rohlenko~\cite{CLR}, \cite{CCL},
have made a very detailed analysis
of some special models with time-periodic potentials.
More precisely, they have found and applied a criterion that
ensures scattering of the wave function.
On the level of the Floquet operator this means that there is no discrete spectrum.
It would be interesting to obtain dispersive estimates for these cases.

Another well-studied class of time-dependent potentials are the so-called charge transfer models.
These are Hamiltonians of the form
\[ H(t) = -\Laplace + \sum_{j=1}^m V_j(\cdot-v_j t) \]
where $\{ v_j\}_{j=1}^m$ are distinct velocities and $V_j$ are
well-localized potentials. They admit localized states that travel
with each of these potentials and asymptotically behave like the
sum of bound states of each of the ``channel Hamiltonians''
 \[ H(t) = -\Laplace +  V_j(\cdot-v_j t)\ . \]
Those are of course Galilei transformed bound states of the corresponding stationary Hamiltonians.
Yajima~\cite{Y} and Graf~\cite{Gr} proved
that these Hamiltonians are asymptotically complete, i.e., that as $t\to\infty$ each state decomposes
into a sum of wave functions associated with each of the channels, including the free channel.

Rodnianski, Soffer, and the author obtained dispersive estimates for these models in the spaces
$L^1\cap L^2\to L^2\cap L^\infty$.
Later, Cai~\cite{cai1} as part of his Caltech thesis removed $L^2$ from these bounds.
Such estimates were needed in order to prove
asymptotic stability of $N$-soliton solutions, see~\cite{RSS2}.

\bibliographystyle{amsplain}

\end{document}